\documentclass[10pt]{article}

\newif\ifdebug                                                      %
\debugtrue

\usepackage{amsmath,amssymb,amsthm}
\usepackage{rotating,pdflscape,afterpage}
\usepackage[sc,osf]{mathpazo}
\AtBeginDocument{
	\DeclareSymbolFont{AMSb}{U}{msb}{m}{n}
	\DeclareSymbolFontAlphabet{\mathbb}{AMSb}
}
\usepackage{graphicx}
\usepackage{extpfeil}
\usepackage{mathabx} 
\usepackage{fancyvrb}
\usepackage{float}
\usepackage{caption}
\newcommand{\mockalph}[1]{\!}
\usepackage{enumitem}
\usepackage[nouppercase]{scrlayer-scrpage}
\usepackage{avant} 
\usepackage{mathrsfs}
\usepackage{mathtools}
\usepackage{faktor,xfrac}
\usepackage[all]{xypic}\xyoption{rotate}
\newdir{ >}{{}*!/-7pt/\dir{>}}
\newdir{i}{{}*!/-7pt/\dir{(}	}

\usepackage{thmtools}

\usepackage[utf8]{inputenc}
\usepackage[T1]{fontenc}
\usepackage{tocloft}
\makeatletter
\renewcommand{\l@figure}{\@dottedtocline{1}{1em}{3.5em}}
\renewcommand{\l@table}{\@dottedtocline{2}{1em}{3.5em}}
\makeatother
\newcommand*{\noaddvspace}{\renewcommand*{\addvspace}[1]{}}
\addtocontents{lof}{\protect\noaddvspace}

\usepackage[hyphens]{url}		%
\usepackage{hyperref}	
\usepackage{cleveref} 	

\makeatletter
\let\c@figure\c@table
\let\c@equation\c@table
\makeatother

\numberwithin{table}{section}
\numberwithin{equation}{section}

\newtheorem{theorem}[table]{Theorem}

\newtheorem{proposition}[table]{Proposition}

\newtheorem{corollary}[table]{Corollary}

\newtheorem{lemma}[table]{Lemma}

\newtheorem{claim}[table]{Claim}

\theoremstyle{definition}

\newtheorem{definition}[table]{Definition}

\newtheorem{notation}[table]{Notation}

\newtheorem{observation}[table]{Observation}

\newtheorem{conjecture}[table]{Conjecture}

\theoremstyle{remark}
\newtheorem{fact}[table]{Fact}

\newtheorem{example}[table]{Example}
\newtheorem{exercise}[table]{Exercise}

\newtheorem{problem}[table]{Problem}
\newtheorem{histrmks}[table]{Historical remarks}
\newtheorem{remark}[table]{Remark}
\newtheorem{remarks}[table]{Remarks}

\theoremstyle{plain}
\newtheorem*{thm*}{Theorem}
\newtheorem*{theorem*}{Theorem}
\newtheorem*{prop*}{Proposition}
\newtheorem*{proposition*}{Proposition}
\newtheorem*{lemma*}{Lemma}
\newtheorem*{corollary*}{Corollary}
\newtheorem*{cor*}{Corollary}

\theoremstyle{definition}
\newtheorem*{definition*}{Definition}
\newtheorem*{defn*}{Definition}
\newtheorem*{QQ*}{Question}
\newtheorem*{obs*}{Observation}
\newtheorem*{notation*}{Notation}

\theoremstyle{remark}
\newtheorem*{rmk*}{Remark}
\newtheorem*{remark*}{Remark}

\newtheorem*{examples*}{Examples}
\newtheorem*{example*}{Example}
\newtheorem*{EG*}{Example}
\newtheorem*{EGs*}{Examples}
\newtheorem*{fact*}{Fact}
\newtheorem*{prob*}{Problem}

\usepackage{etoolbox}
\usepackage[usenames,dvipsnames,svgnames,table]{xcolor}
\newcommand		{\defd}[1]	{\textcolor{RoyalBlue}{\textbf{\textit{#1}}}}
\newcommand		{\defm}[1]	{\textcolor{RoyalBlue}{#1}}

\tracingpatches
\makeatletter
\patchcmd{\@setref}{\bfseries ??}{\bfseries\color{red} FIX ME!}{}{}
\patchcmd{\@setcite}{\bfseries ?}{\bfseries\color{red} FIX ME!}{}{}
\patchcmd{\@setcref}         {??}{\color{red} FIX ME!}{}{}
\patchcmd{\@setcref}         {??}{\color{red} FIX ME!}{}{}
\patchcmd{\@setcrefrange}    {??}{\color{red} FIX ME!}{}{}
\patchcmd{\@setcrefrange}    {??}{\color{red} FIX ME!}{}{}
\patchcmd{\@setcrefrange}    {??}{\color{red} FIX ME!}{}{}
\patchcmd{\@setcrefrange}    {??}{\color{red} FIX ME!}{}{}
\patchcmd{\@setcrefrange}    {??}{\color{red} FIX ME!}{}{}
\patchcmd{\@setcrefrange}    {??}{\color{red} FIX ME!}{}{}
\patchcmd{\@setnamecref}     {??}{\color{red} FIX ME!}{}{}
\patchcmd{\@setnamecref}     {??}{\color{red} FIX ME!}{}{}
\patchcmd{\@setcpageref}     {??}{\color{red} FIX ME!}{}{}
\patchcmd{\@setcpageref}     {??}{\color{red} FIX ME!}{}{}
\patchcmd{\@setcpagerefrange}{??}{\color{red} FIX ME!}{}{}
\patchcmd{\@setcpagerefrange}{??}{\color{red} FIX ME!}{}{}
\patchcmd{\@setcpagerefrange}{??}{\color{red} FIX ME!}{}{}
\patchcmd{\@setcpagerefrange}{??}{\color{red} FIX ME!}{}{}
\patchcmd{\@setcpagerefrange}{??}{\color{red} FIX ME!}{}{}
\patchcmd{\@cref}            {??}{\color{red} FIX ME!}{}{}
\def\blx@citation@entry#1#2{%
	\blx@bibreq{#1}%
	\ifinlist{#1}{\blx@cites}
	{}
	{\listgadd{\blx@cites}{#1}%
		\blx@auxwrite\@mainaux{}{\string\abx@aux@cite{#1}}}%
	\ifinlistcs{#1}{blx@segm@\the\c@refsection @\the\c@refsegment}
	{}
	{\listcsgadd{blx@segm@\the\c@refsection @\the\c@refsegment}{#1}}%
	\blx@ifdata{#1}%
	{}%
	{\ifcsdef{blx@miss@\the\c@refsection}%
		{\ifinlistcs{#1}{blx@miss@\the\c@refsection}%
			{{\bfseries\color{red} cite:} }%
			{\blx@logreq@active{#2{#1}}}}%
		{\blx@logreq@active{#2{#1}}}}}
\def\blx@citeadd#1{%
	\ifcsdef{blx@keyalias@\the\c@refsection @#1}
	{\edef\blx@realkey{\csuse{blx@keyalias@\the\c@refsection @#1}}}
	{\def\blx@realkey{#1}}%
	\expandafter\blx@citation\expandafter{\blx@realkey}\blx@msg@cundefon
	\expandafter\blx@ifdata\expandafter{\blx@realkey}
	{\advance\blx@tempcnta\@ne
		\listeadd\blx@tempa{\blx@realkey}}
	{\ifnum\blx@tempcntb>\z@MUlticitedelim\fi
		\expandafter\abx@missing\expandafter{\blx@realkey}%
		\advance\blx@tempcntb\@ne}}

\DeclarePairedDelimiterX{\pmodx}[1]{(}{)}{{\operator@font mod}\mkern6mu#1}
\renewcommand{\pmod}{%
	\allowbreak
	\if@display\mkern18mu\else\mkern8mu\fi
	\pmodx
}
\makeatother

\makeatletter
\newcommand{\oset}[3][0ex]{%
	\raisebox{.175ex}{$%
		\mathrel{\mathop{#3}\limits^{
				\vbox to#1{\kern-2\ex@
					\hbox{$\scriptstyle#2$}\vss}}}
		$}%
}
\makeatother

\newcommand{\myred}{BrickRed}
\hypersetup{
	colorlinks			= true, 	
	urlcolor			= \myred, 
	linkcolor			= \myred, 	
	citecolor			= PineGreen 	
}

\usepackage{xspace}

\usepackage[left=2.54cm,right=2.54cm,top=2.54cm,bottom=2.54cm]{geometry}
\usepackage{tikz}
\usetikzlibrary{matrix,fit,backgrounds,calc,arrows} 

\tikzstyle{image}=[rectangle,fill=Red!20,inner sep=-2pt]
\tikzstyle{nonzero}=[rectangle,fill=Navy!20,inner sep=0pt]
\tikzstyle{nonzerosm}=[rectangle,fill=Navy!20,inner sep=-2pt]

\makeatletter
\newbox\xrat@below
\newbox\xrat@above
\newcommand{\xrightarrowtail}[2][]{%
	\setbox\xrat@below=\hbox{\ensuremath{\scriptstyle #1}}%
	\setbox\xrat@above=\hbox{\ensuremath{\scriptstyle #2}}%
	\pgfmathsetlengthmacro{\xrat@len}{max(\wd\xrat@below,\wd\xrat@above)+.6em}%
	\mathrel{\tikz [>->,baseline=-.55ex]
		\draw (0,0) -- node[below=-2pt] {\box\xrat@below}
		node[above=-2pt] {\box\xrat@above}
		(\xrat@len,0) ;}}
\newbox\xrat@below
\newbox\xrat@above
\renewcommand{\xtwoheadrightarrow}[2][]{%
	\setbox\xrat@below=\hbox{\ensuremath{\scriptstyle #1}}%
	\setbox\xrat@above=\hbox{\ensuremath{\scriptstyle #2}}%
	\pgfmathsetlengthmacro{\xrat@len}{max(\wd\xrat@below,\wd\xrat@above)+.6em}%
	\mathrel{\tikz [->>,baseline=-.55ex]
		\draw (0,0) -- node[below=-2pt] {\box\xrat@below}
		node[above=-2pt] {\box\xrat@above}
		(\xrat@len,0) ;}}
\makeatother

\newcommand{\xepi}{\xtwoheadrightarrow}

\makeatletter
\newcommand{\presectionskip}{-1.5\baselineskip}
\newcommand{\postsectionskip}{0.3\baselineskip}
\usepackage{titlesec}
\renewcommand{\section}{\@startsection
	{chapter}{0}{0mm}
	{\presectionskip}
	{\postsectionskip}
	{\sffamily\huge}}
\renewcommand{\section}{\@startsection
	{section}{1}{0mm}
	{\presectionskip}
	{\postsectionskip}
	{\sffamily\LARGE}}
\renewcommand{\subsection}{\@startsection
	{subsection}{2}{0mm}
	{\presectionskip}
	{\postsectionskip}
	{\sffamily\Large}}
\renewcommand{\subsubsection}{\@startsection
	{subsubsection}{3}{0mm}
	{\presectionskip}
	{\postsectionskip}
	{\sffamily\normalsize}}
\renewcommand{\@seccntformat}[1]{\csname the#1\endcsname.\quad}
\newcommand\HUGE{\@setfontsize\Huge{30}{47}} 
\titleformat{\chapter}[display]
{\sffamily\Large}
{Chapter {\HUGE\normalfont\thechapter}}    
{1em}
{\huge}
\titleformat{\section}
{\sffamily\Large}
{\normalfont{\@setfontsize\Large{18}{47}\thesection.}}    
{1em}
{}
\titleformat{\subsection}
{\sffamily\large}
{\normalfont{\Large\thesubsection}.}    
{1em}
{}
\titleformat{\subsubsection}
{\sffamily\normalsize}
{\normalfont{\large\thesubsubsection}.}    
{1em}
{}
\makeatother

\makeatletter
\def\smallunderbrace#1{\mathop{\vtop{\m@th\ialign{##\crcr
				$\hfil\displaystyle{#1}\hfil$\crcr
				\noalign{\kern3\p@\nointerlineskip}%
				\tiny\upbracefill\crcr\noalign{\kern3\p@}}}}\limits}
\makeatother

\renewcommand{\SS}{\textsection}
\newcommand{\bthm}{\begin{theorem}}
	\newcommand{\ethm}{\end{theorem}}
\newcommand{\bprop}{\begin{proposition}}
	\newcommand{\eprop}{\end{proposition}}
\newcommand{\bcor}{\begin{corollary}}
	\newcommand{\ecor}{\end{corollary}}
\newcommand{\bconj}{\begin{conjecture}}
	\newcommand{\econj}{\end{conjecture}}
\newcommand{\blem}{\begin{lemma}}
	\newcommand{\elem}{\end{lemma}}
\newcommand{\bclm}{\begin{claim}}
	\newcommand{\eclm}{\end{claim}}
\newcommand{\bpf}{\begin{proof}}
	\newcommand{\epf}{\end{proof}}
\newcommand{\bdetails}{\begin{details}}
	\newcommand{\edetails}{\end{details}}
\newcommand{\bdefi}{\begin{definition}}
	\newcommand{\edefi}{\end{definition}}
\newcommand{\bdefn}{\begin{definition}}
	\newcommand{\edefn}{\end{definition}}
\newcommand{\bex}{\begin{example}}
	\newcommand{\eex}{\end{example}}
\newcommand{\bprob}{\begin{problem}}
	\newcommand{\eprob}{\end{problem}}
\newcommand{\bob}{\begin{observation}}
	\newcommand{\eob}{\end{observation}}
\newcommand{\bexer}{\begin{exercise}}
	\newcommand{\eexer}{\end{exercise}}
\newcommand{\bexers}{\begin{exercises}}
	\newcommand{\eexers}{\end{exercises}}
\newcommand{\brmk}{\begin{remark}}
	\newcommand{\ermk}{\end{remark}}
\newcommand{\bhist}{\begin{histrmks}}
	\newcommand{\ehist}{\end{histrmks}}
\newcommand{\brmks}{\begin{remarks}}
	\newcommand{\ermks}{\end{remarks}}
\newcommand{\bverb}{\begin{verbatim}}
	\newcommand{\everb}{\end{verbatim}}
\newcommand{\bntn}{\begin{notation}}
	\newcommand{\entn}{\end{notation}}
\newcommand{\bfct}{\begin{fact}}
	\newcommand{\efct}{\end{fact}}
\newcommand{\bfcts}{\begin{facts}}
	\newcommand{\befcts}{\end{facts}}
\newcommand{\benum}{\begin{enumerate}}
	\newcommand{\eenum}{\end{enumerate}}
\newcommand{\bitem}{\begin{itemize}}
	\newcommand{\eitem}{\end{itemize}}

\renewcommand	{\o}		{\circ}

\renewcommand	{\epsilon}	{\varepsilon}

\newcommand		{\h}		{\eta}

\renewcommand	{\:}		{\colon}

\newcommand		{\IpG}	{J}
\newcommand		{\oIpG}	{\oplus\IpG}

\newcommand		{\AHSS}	{Atiyah--Hirzebruch spectral sequence\xspace}

\renewcommand	{\th}			{^{\mathrm{th}}}

\newcommand{\compl}{\!\wh{\ ^{\phantom{x}}}\mn}
\newcommand{\K}{K^*}
\newcommand{\KG}{\K_G}

\makeatletter
\newcommand{\subalign}[1]{%
	\vcenter{%
		\Let@ \restore@math@cr \default@tag
		\baselineskip\fontdimen10 \scriptfont\tw@
		\advance\baselineskip\fontdimen12 \scriptfont\tw@
		\lineskip\thr@@\fontdimen8 \scriptfont\thr@@
		\lineskiplimit\lineskip
		\ialign{\hfil$\m@th\scriptstyle##$&$\m@th\scriptstyle{}##$\crcr
			#1\crcr
		}%
	}
}
\makeatother

\newcommand		{\eqn}[1]			{\begin{align*} #1 \end{align*}}
\newcommand		{\quation}[1]		{\begin{equation} #1 \end{equation}}

\newcommand		{\case}[1]			{\begin{cases} #1 \end{cases}}

\newcommand		{\hyref}[1]			{\hyperref[#1]{\ref*{#1}}} 

\newcommand		{\bs}				{\bigskip}

\newcommand		{\mn}				{\mspace{-2mu}}

\newcommand		{\dsp}				{\displaystyle}

\newcommand		{\nd}				{\noindent}
\newcommand		{\ol}				{\overline}
\newcommand		{\os}			{\overset}
\newcommand		{\us}			{\underset}

\newcommand		{\ul}			{\underline}
\newcommand		{\wh}			{\widehat}

\newcommand		{\wt}			{\widetilde}

\newcommand		{\mr}			{\mathrm}
\newcommand		{\bb}			{\mathbb}
\newcommand		{\mb}			{\mathbf}

\newcommand		{\ms}			{\mathscr}

\newcommand		{\g}			{\gamma}
\newcommand		{\e}			{\epsilon}

\newcommand		{\s}		{\sigma}

\newcommand		{\W}		{\Omega}

\DeclareSymbolFont{cmletters}{OT1}{cmr}{m}{n}
\DeclareMathSymbol{\Ups}{\mathalpha}{cmletters}{"7}
\renewcommand	{\Upsilon}{\Ups}

\newcommand		{\0}		{\varnothing}

\newcommand		{\F}		{\bb F}

\newcommand		{\Z}		{\bb Z}
\newcommand		{\Q}		{\bb Q}
\newcommand		{\R}		{\bb R}
\newcommand		{\C}		{\bb C}

\newcommand		{\RP}		{\bb R \mr P}
\newcommand		{\CP}		{\bb C \mr P}
\newcommand		{\CPi}		{\CP^\infty}

\DeclareMathOperator{\id}	{id}

\renewcommand 	{\H}	{H^*}
\newcommand 	{\HG}	{{\H_G}}

\newcommand{\MU}{MU^*}

\let\inter\cap%

\newcommand		{\dis}		{\amalg}

\newcommand		{\6}			{\partial}

\renewcommand	{\-}		{^{-1}}
\renewcommand	{\o}		{\circ}
\renewcommand	{\.}		{\cdot}
\newcommand		{\x}		{\times}

\newcommand{\oplushigher}{\mathbin{\raisebox{.85pt}{$\displaystyle\oplus$}}}

\DeclareMathOperator*{\otimesvariable}{%
	\mathchoice {\raisebox{.85pt}{$\displaystyle\otimes$}}
	{\raisebox{.85pt}{$\otimes$}}
	{\raisebox{0.7pt}{$\scriptstyle\otimes$}}
	{\raisebox{0.2pt}{$\scriptscriptstyle\otimes$}}
}
\newcommand		{\tensor}		{\otimesvariable}
\newcommand		{\ox}			{\tensor}

\newcommand		{\direct}		{\oplushigher}

\newcommand		{\+}			{\direct}

\newcommand		{\Direct}		{\bigoplus}

\newcommand		{\colim}		{\varinjlim}

\DeclareMathOperator{\Hom}		{Hom}

\newcommand		{\U}			{\mr{U}}

\newcommand{\Th}		{\mr{Th}}

\newcommand		{\longto} 		{\longrightarrow}
\newcommand		{\lt}			{\longto}

\newcommand		{\lmt}			{\longmapsto}

\newcommand		{\inc}		{\hookrightarrow}
\newcommand		{\xinc}		{\xhookrightarrow}
\newcommand		{\longinc}		{\xinc[]{\ \ \ \ }}

\newcommand		{\longepi}	{\xepi[]{\ \ \ \ }}

\newcommand		{\simto}		{\xrightarrow{\sim}}

\newcommand		{\longsimto}	{\os\sim\longto}
\newcommand		{\isoto}		{\longsimto}

\newcommand		{\ceq}			{\coloneqq}
\newcommand		{\eqc}			{\eqqcolon}

\newcommand		{\hmt}			{\simeq}
\newcommand		{\iso}				{\cong}

\usepackage{todonotes}

\theoremstyle{definition}

\newcommand{\AID}{abstract isotropy data\xspace}

\renewcommand{\sp}{\s_{\mn p}}

\newcommand{\AP}{[\ms A\!, \ms P]}
\newcommand{\mfring}{\Omega_*^{U:T}}
\newcommand{\bundlering}{\Omega_*^{U:T}[\ms A\!,\ms P]}
\newcommand{\WUG}{\Omega^{U:G}_*}
\newcommand{\MUG}{MU^G_*}
\newcommand{\WUT}{\Omega^{U:T}_*}
\newcommand{\MUT}{MU^T_*}
\newcommand{\WUTAP}{\Omega^{U:T}_*[\ms A\!, \ms P]}
\newcommand{\WUGAP}{\Omega^{U:G}_*[\ms A\!, \ms P]}
\newcommand{\MUTAP}{MU^T_*[\ms A\!, \ms P]}
\newcommand{\MUGAP}{MU^G_*\AP}
\newcommand{\WSF}{\Omega_*^{\mathrm{SF}}}
\newcommand{\MUSF}{MU_*^{\mathrm{SF}}}

\begin{document}
	
	\title{\vspace{-1em}\huge Fixed points
		and bordism of 
		semifree actions}
	\author{Jeffrey D.~Carlson}
	\maketitle

	\begin{abstract}%
		\small{%
			We apply fixed-point techniques to compute the coefficient ring
			of semifree circle-equivariant complex cobordism with isolated fixed points, 
			recovering a 2004 result of Sinha through 19$^{\textrm{th}}$-century methods. 
			
			%
		}
	\end{abstract}

	Homotopical equivariant complex cobordism $MU^{G}$ 
	with respect to the action of a compact Lie group
	$G$ is the universal complex-oriented $G$-equivariant spectrum,
	yet concrete presentations of its coefficient ring are 
	unknown except for certain cases when $G$ is finite.
	Geometric equivariant cobordism $\Omega^{U:G}_*(-)$, 
	whose coefficient ring
	is given by bordism classes of (tangentially) stably complex closed $G$-manifolds
	is similarly inaccessible in most cases
	as of this writing.

	Restricting to the special case 
	of \emph{semifree} circle actions,
	where all orbits are free or fixed points, 
	Sinha~\cite{sinha2005semifree} proved a number of results,
	including, strikingly, the following.
	\begin{theorem}[Sinha~{\cite[Thm.~1.1]{sinha2005semifree}}]\label{thm:semifree bordism}
		Every compact, oriented, stably complex, semifree $S^1$-manifold 
		with isolated fixed points
		is cobordant to a disjoint union of direct powers of 
		the sphere $S^2 = \CP^1$ 
		with the standard complex structure 
		and rotation action.
		That is, the bordism ring  
		of such manifolds is isomorphic to the polynomial ring $\Z[S^2]$ on one generator.
	\end{theorem}
\nd	In the present note, we recover this result 
	through elementary means.

	\section{Introduction}
We begin by reviewing the definitions and existing results.\footnote{\ 
	The literature review in the first version of this document 
	was deliberately sketchy, with a view toward brevity,
	and several correspondents requested it be expanded and clarified.
	The corresponding section of a second version grew to a survey of twenty-some pages, which remains unfinished.
	One came to realize it should be separate, 
	if it is ever finished at all; 
	there is clearly a balance to be struck in a note of this length. 
	Let us try again.
	}
We will not attempt to discuss the role of equivariant
complex cobordism in equivariant homotopy theory, 
but will try to state, briefly, what is known about it
from a computational perspective.
	We take as a starting point Thom's result that 
	(nonequivariant) geometric complex bordism $\defm{\Omega^U_n}$, 
	defined in terms of bordism of stably complex $n$-manifolds,
 	is also represented by the 
 	Thom spectrum $\defm {MU}$
	whose $2n\th$ level $\defm {MU_{2n}}$ is the Thom space
	of the universal complex $n$-plane bundle, 
	with the isomorphism induced from the Pontrjagin--Thom collapse.\footnote{\ 
		We also assume
	the classic computational result, 
	due to Milnor and Novikov,
	that 
	$\defm{\Omega^U_*} \iso \defm{MU_*} \ceq \pi_* MU$  
	is a polynomial ring
	$\Z[x_2,x_4,x_6,\ldots]$ on one generator $x_{2d} \in \pi_{2d} MU$ 
	for each positive natural $n$.
	One has $\Omega^*_U \ox \Q \iso \Q\big[[\CP^1],[\CP^2],[\CP^3],\ldots\big]$, 
	though not in such a way that the $[\CP^d]$ are scalar multiples of 
	the integral generators $x_{2d}$;
	Milnor found a natural set of manifolds giving integral generators,
	but the relations they satisfy are not straightforward.
	The famous theoretical result, due to Quillen, 
	is that this coefficient ring $MU_*$ 
	carries Lazard's universal one-dimensional
	formal group law. 
	It is also important, but more tautological,
	that $MU^*(-)$ is initial amongst
	complex-oriented cohomology theories.
}
	The corresponding map for equivariant complex bordism is not an 
	isomorphism, leading to distinct geometric and homotopical notions.
	
\subsection{Definitions and comparison}
	
	Fix a compact Lie group $\defm G$.
	A $G$-equivariant complex vector bundle (or \emph{$G$-bundle})
	over a $G$-space $M$
	is the usual thing, a complex vector bundle $E \to M$
	equipped with a $\C$-linear $G$-action on $E$
	such that the projection is equivariant.
	We write $\defm{\R^k}$ (resp. $\defm{\C^k}$)
	for Euclidean space viewed as a trivial real 
	(resp. complex) 
	$G$-representation,
	and when $V$ is a $G$-representation
	and we are considering some $G$-manifold $M$,
	we denote by $\defm{\ul V}$ the trivial bundle $M \x V \to M$,
	equipped with the diagonal $G$-action on its total space. 
	By a \defd{stably complex $G$-manifold} we mean a smooth manifold $M$
	together with a complex $G$-bundle structure
	on some stabilization $TM \+ \ul\R^k$ of its tangent $G$-bundle.
	Two $G$-equivariant stable complex structures on a manifold are defined to be \emph{equivalent}
	if after direct-summing $\ul\C^m$ to one and $\ul\C^n$
	to the other, for some $m$ and $n$, 
	the two become isomorphic as complex $G$-bundles 
	via an isomorphism fixing $TM$.
		More concisely, an equivalence class 
		is a lift of the class of $[TM]$
		in reduced real equivariant topological K-theory 
		${KO}{}^0_G (M)/\Z$
		along the forgetful map from reduced complex equivariant K-theory
		$K^0_G (M) /\Z$,
	where the two instances of $\Z$ are the summands 
	respectively generated by $\ul\R$ and $\ul\C$.
The \defd{equivariant geometric bordism} groups $\defm{\Omega^{U:G}_{n}}(X)$ 
are determined 
as singular $G$-maps of stably complex $G$-manifolds $M^n$ into $X$,
modulo the equivalence relation identifying two such maps
if they together form the boundary restriction of a $G$-map of a stably complex $G$-manifold $W^{n+1}$
into $X$.

Different notions of lifting and equivalence yield different groups:
one might instead ask of $TM$ that
there exist a real representation
$W$ such that $TM \oplus \ul W$ is a complex $G$-bundle
and define two such structures 
$TM \+ \ul W$ and $TM \+ \ul W'$ to be equivalent
if there exist complex $G$-representations $\ul V,\ul V'$
such that 
$TM \+ \ul W \+ \ul V$ and $TM \+ \ul W' \+ \ul V'$
are isomorphic as complex $G$-bundles.
	In K-theoretic terms, an element is then given by a lift of
	$[TM] \in \smash{\wt{KO}{}^0_G} (M)$
	in ${\wt K}^0_G (M)$.		
Equivalently, 
one can ask the same of the normal bundle 
$\nu =s{\nu_V(M)}$ to $M$
under an equivariant embedding 
of $M$ in a sufficiently large real $G$-representation $V$;
but as $[\nu] = -[TM]$ in $\wt{K}^0_G(M)$,
these are equivalent notions. 
The resulting groups have variously been denoted 
$\defm{mu^{G}_*}(X)$ 
and $\defm{\ol\W{}^{U:G}_*}(X)$;
see Darby's thesis~\cite[\SS\SS3.1--2]{darbythesis}
and Comeza\~na~\cite[\SS3]{comezana}.\footnote{\ 
	An analogous discussion can be made of 
	equivariant	unoriented bordism,
	spawning theories $\W^{O:G}_*$, $mo^G_*$, and $MO^G_*$
	which very frequently admit analogous results.
	Indeed, many proofs carry over with little to no change
	from one setting to the other.
	Perhaps heretically,
	but with the aim of bounding our exposition,
	we bypass any discussion of these theories.
}

\bex[{Hanke~\cite[p.680]{hankebordism}}]
To see these notions are distinct, consider the 
three-dimensional representation $V$ of $\Z/2$
determined by the antipodal map.
Its unit sphere $\defm{SV}$ is invariant and hence a $\Z/2$-space,
and the normal bundle $\nu_{V}(SV)$ is the trivial $\Z/2$-bundle $\ul \R$
over $SV$, which can be stabilized to $\ul \C$.
On the other hand, the tangent $\Z/2$-bundle $TSV$
is naturally a pullback of the tangent bundle $T\RP^2$,
so a complex $\Z/2$-bundle structure on $TSV \+ \ul\R^{2k}$
would amount to a complex structure on $T\RP^2 \+ \ul \R^{2k}$.
But that in turn would induce a complex structure
on $T(\RP^2 \x \R^{2k})$, which is impossible
as $\RP^2 \x \R^{2k}$ is not orientable.
\eex

There is at least a natural map 
$\WUG(X) \lt mu^{G}_*(X)$,
which L\"offler asserts~\cite[Lem.~2.1]{loefflercharacteristic}
and a result of Comeza\~na~\cite[Thm.~5.4]{comezana}
implies is injective for all compact abelian $G$.
Writing $\defm{S^V}$ for the one-point compactification
of a complex $G$-representation $V$
and $\defm{|V|}$ for its real dimension,
one can stabilize both $\Omega^{U:G}_*$
and $mu^G_*$ 
by maps
$\smash{\W^G_k(X) \lt \W^G_{k+|V|}(S^V \wedge X)}$	
and 
$\smash{mu^G_k(X) \lt mu^G_{k+|V|}(S^V\wedge X)}$
given in terms of a representative $M \to X$
by ``straightening the angles'' of $S^V \wedge M$
to make it again a manifold, 
and it is a result of Comeza\~na--Costenoble~\cite[Thm.~3.3]{comezana} 
that the colimit of the maps 
$\smash{
	\Omega^{U:G}_{k+|V|}(S^V \wedge X)
} 
	\lt 
\smash{
	mu^G_{k+|V|}(S^V \wedge X)
}$
is an isomorphism.\footnote{\ 
	N.B. the typo on p.~392: 
	the domain and codomain
	of $\phi$ should be transposed
	and similarly for $\Phi$.
	}
This colimit gives a third variant of equivariant complex bordism.

We claim this homology theory is representable,
unlike $\WUG(-)$ and $mu^G_*(-)$.
To see this, note that there exists a \emph{universal}
complex $G$--vector bundle
$\defm{EU_{2n}^G \to BU_{2n}^G}$
such that homotopy classes of $G$-maps $X \lt BU^G_{2n}$
parameterize isomorphism classes of
rank-$n$ complex $G$-bundles over a $G$-space $X$;
the most natural model is the tautological bundle over
the Grassmannian of complex 
$n$-planes in the direct sum 
$\defm{\mathbf{U}} = \Direct V^\infty$ of $\aleph_0$-many 
copies of each irreducible complex $G$-representation $V$.
The associated Thom spaces $\defm{MU^G_{2n}}$ 
admit
a multiplication $MU^G_{2k} \wedge MU^G_{2m} \lt MU^G_{2k+2m}$ 
defined from the maps
$BU^G_{2k} \x BU^G_{2m} \lt BU^G_{2k+2m}$
classifying the Whitney sum.
For any representation $V$ of $G$, the classifying map
$\smash{\ul V \+ EU^G_{2k} \lt EU^G_{|V|+2k}}$ induces a map of Thom spaces
$\smash{S^V \wedge MU^G_{2k}} \lt \smash{MU^{\smash G}_{|V|+2k}}$
which on forgetting the $G$-action
becomes the structure map for the nonequivariant
sequential Thom spectrum $MU$. 
We take these as structure maps for the $G$-spectrum $\defm{MU^G}$.\footnote{\ 
	One can expand this to an $R(G)$-graded spectrum
	in the Lewis--May--Steinberger formalism
	without much trouble
	by taking $\defm{BU_V^G}$ 
	to be the Grassmannian of complex $\frac 1 2|V|$-planes
	in the abstract vector space $V \+ \mathbf{U}$,
	basepointed at $V$,
	and one gets $\smash{MU_{V}^G(-)} \iso \smash{MU_{|V|}^G(-)}$
	and similarly for cobordism.
	The structure maps 
	$\smash{S^{V/W} \wedge MU^G_W \lt MU^G_V}$
	for $W < V$ arise from applying the Thom construction
	to the map of tautological bundles
	taking each complex $\frac 1 2|W|$-plane $L$ in $W \+ \mb U$
	to the $\frac 1 2|V|$-plane ${V/W} \+ L$ in $V \+ \mb U$,
	where $V/W$ here denotes the orthogonal complement
	with respect to the inner product on $\mb U$.

	It is asserted semifrequently in the literature that 
	one can expand this to an $RO(G)$-graded spectrum,
	but details never seem to follow.
	Here are some. 
	We index by \emph{real} $G$-invariant subspaces 
	of the complex $G$-universe $\mb U$.
	For each such space $W$, write $W^\C \ceq W \inter iW$,
	which is a canonical maximal complex subrepresentation,
	and set
	$\smash{MU^G_W \ceq S^{W/W^\C} \wedge MU^G_{W^\C}}$,
	where $W/W^\C$ denotes the orthogonal complement
	and $MU^{G\phantom{'}}_{W^\C}$ is as in the previous paragraph.
Then structure maps
	$S^{V/W} \wedge MU^G_W \lt MU^G_{V}$
	are defined for $W < V$ by 
	first performing the reassociation
$\smash{
	S^{V/W} \wedge S^{W/W^\C} 
		\iso
	S^{V/W^\C} 
	\iso
	S^{V/V^\C} \!\!\wedge S^{V^\C/W^\C}
}\,	$ in the domain,
	then leaving the factor $\smash{S^{V/V^{\C}}}$
	alone and applying the structure map
	$S^{V^\C/W^\C} \wedge  MU^G_{W^\C}
		\lt
	MU^G_{V^\C}$ from the $R(G)$-graded definition
	to the other factor.
	In words, one applies as much of the $R(G)$-graded
	structure as one can, benignly neglecting a 
	``purely real'' remainder sphere.
}
Once $\mathbf{U}$ is equipped with an invariant Hermitian inner product,
each pair $W \leq V \leq \mb U$
defines a unique complementary representation $V/W \leq V$.
We may then partially order $G$-representations 
and use these complements and the structure maps
just defined to associate reduced (co)homology theories 
by
\eqn{\label{def:MUG}
		\defm{\wt{MU}{}_G^{2n}(X)} 
		\ceq 
				\colim_V \,[S^V \wedge X, MU^G_{2n+|V|}]^G\mathrlap,
			\qquad\qquad
		\defm{\wt{MU}^{\smash G}_{2n}(X)} 
		\ceq 
			\colim_V \,[S^V, X \wedge MU^G_{|V|-2n}]^G\mathrlap,
}
where $X$ is a space with a $G$-fixed basepoint
and
$\defm{[-,-]^G}$ denotes pointed $G$-homotopy classes of $G$-maps.
In odd dimensions we instead set
$\defm{\wt{MU}{}^{2n-1}_G(X)} \ceq \wt{MU}{}^{2n}_{\smash G}(S^1 \wedge X)$
and
$\defm{\wt{MU}{}_{2n-1}^G(X)} \ceq \wt{MU}{}_{2n}^{\smash G}(S^1 \wedge X)$.
The unreduced versions 
are defined by
$\defm{MU{}^*_G(X)} \ceq \wt{MU}{}^*_G(X_+)$
and
$\defm{MU{}_*^G(X)} \ceq \wt{MU}{}_*^G(X_+)$
where $X_+$ is the disjoint union of $X$ with a 
new $G$-fixed basepoint.
This is \defd{homotopical equivariant complex (co)bordism}.

Geometric equivariant bordism is related to the homotopical variety via
the Pontrjagin--Thom collapse:
if $M \to X$ represents an element of $mu_n^G(X)$,
then letting $M \longinc V$ be an equivariant embedding
whose normal $G$-bundle $\nu$ carries a complex structure,
the collapse
$S^V \longepi \Th\, \nu$
followed by
the projection $\nu \longepi M$ 
and the classifying map of $\nu$
induce a map $S^V \lt X_+ \wedge MU^G_{|V|-|M|}$,
where $\defm{|M|}$ is the real dimension of $M$.
One checks that increasing $V$ to $V \+ W$
in the definition of the normal embedding
gives the $W$-suspension of the first map,
so that the image in the colimit is well-defined.
The theorem of Br\"ocker--Hook~\cite{broeckerhook}\footnote{\ 
	proven for $mo^G_* \lt MO^G_*$,
	but the proof carries over \emph{mutatis mutandis}
} 
is that stabilizing along these maps induces an isomorphism
\[
\colim_V \,mu^G_{k+|V|}(S^V \wedge X) \isoto\wt{MU}{}^G_{k}(X).\footnote{\ 
	This stabilization is in fact necessary to make equivariant cobordism
	representable, for the suspension maps are not generally isomorphisms
	before passing to the colimit.
}
\]
Composing this isomorphism with the Comeza{\~n}a--Costenoble
isomorphism, one sees $\wt{MU}{}_*^G(-)$
is the stabilization of $\W_*^{U:G}(-)$
as well.
That the map $\defm\Psi\: \W_*^{U:G}(-) \lt \wt{MU}{}_*^G(-)$
is not an isomorphism before stabilizing
can be seen as a result of the failure of 
the Thom transversality argument
in the equivariant setting:	
	generally there is no way to \emph{equivariantly} 
	homotope an arbitrary $G$-map to be transversal. 
	
	\bex[{Wasserman~\cite[p.~137]{wassermanequivariant}}]
	Consider $\R$ as a $\Z/2$-manifold under $x \lmt -x$ 
	and the smooth equivariant map $f$ sending $M = {*}$ 
    to the $\Z/2$-submanifold $W = \{0\}$ of $\R$.
	There is no way to equivariantly homotope $f$ to be transverse to $W$;
	in fact there are no other $\Z/2$-homotopic maps at all.
\eex
	
\nd On the other hand,
when the $G$-action on $X$ is free,
the transversality argument goes through and 
$\Omega^{U:G}_*(X) \lt MU^G_*(X)$ 
is an isomorphism~\cite[Prop.~1.3]{tomdieckintegrality}.

\subsection{Comparisons and general structure}\label{sec:comparison}

None of the theories $X \lmt \W^{U:G}_*(X)$, 
$mu^{G}_*(X)$, and
$MU^{G}_*(X)$ is readily calculable, 
and each is only known for special values of $G$ and $X$,
but there are a number of comparison maps with other theories,
very frequently injective, and in many cases identifying
cobordism as a pullback.
Some general structure results are also known.

\subsubsection*{Formal group laws}
Classically, it is known that $MU_*$ carries the universal one-dimensional
formal group law.
Greenlees~\cite{greenleesFGL} 
showed that for $A$ an abelian group,
the classifying map $L_A \lt MU^A_*$
from the ring carrying the universal $A$-equivariant formal group law
is surjective.
Hanke--Wiemeler~\cite{hankewiemeler} 
showed this map is an isomorphism for $A = \Z/2$
and Hausmann~\cite{hausmannglobal} 
showed it for all compact abelian Lie groups $A$.

\subsubsection*{Module structure and evenness}
The map from $MU_*$
assigning a stably complex manifold $M$ the trivial $G$-action
is a ring map inducing an $MU_*$-module structure on
$\Omega^{U:G}_*$, $mu^G_*$, and $MU^*_G$,
and splits the augmentation map to $MU_*$ 
induced by forgetting the action on a $G$-manifold.

Hamrick and Ossa~\cite{hamrickossa} showed that
if $G$ is a topologically cyclic compact Lie group, 
\emph{i.e.}, a product of a finite cyclic group and a torus,
$\Omega^{U:G}_*$ is a free $MU_*$-module on even-degree generators,
after earlier work by several authors establishing special cases.
L\"offler~\cite{loefflercharacteristic} announced
and Comeza\~{n}a~\cite[\SS5]{comezana} 
showed the same of 
$\Omega^{U:G}_*$
and $MU^*_G$ when $G$ is a compact abelian Lie group.
He showed moreover
that if $X$ is a $G$-space of the form
$S^W \x \prod_{j=1}^\ell BU^G_{m_j}$,
then $\W^{U:G}_*(X)$ is free on even-dimensional generators
and
the stabilization maps 
\[
\W^{U:G}_*(BU_n^G \x X) \to 
\W^{U:G}_*(BU_{n+1}^G \x X)\mathrlap,\qquad \ 
\W^{U:G}_*(X) \to 
\W^{U:G}_{*+2}(S^V \x X) \ (V\mbox{ irreducible})\mathrlap,\qquad \ 
\Omega^{U:G}_* \lt MU^*_G
\]
are all split injections of $MU_*$-modules.
The same
is conjectured to hold for all compact Lie groups~\cite{uribeevenness},
and 
Comeza\~{n}a (p.~398) stated he could also prove it 
for the finite dihedral groups and $O(2)$,
but these results appear not to have been published.

\subsubsection*{Reduction to K-theory}

Okonek~\cite{okonekconnerfloyd} showed the natural transformation	
$MU^*_G(-)\lt K^*_G(-)$,
induces a natural isomorphism 
\[
	R(G) \ox_{MU^*_G} MU^*_G(X) \iso \KG (X)
\]
of $\Z/2$-graded rings,
where the module structure map $MU^*_G \lt R(G)$
is the case $X = {*}$.
This generalizes
the Conner--Floyd isomorphism $\Z \ox_{MU^*}   MU^*(X)  \iso \K(X)$%
~\cite{connerfloydK}. 
There are many other such transformations 
$\MUG(-) \lt h_G^*(-)$
(not all so easily characterized),
since $\MUG(-)$ is the universal $G$-equivariant
complex-oriented cobordism theory~\cite{okonekconnerfloyd}.

\subsubsection*{Injections in larger theories}

tom Dieck~\cite{tomdieckintegrality}\footnote{\ 
	His prototype
	was Boardman's map
	$\Omega^{O:\Z/2}_* \lt MO^*(B\Z/2)$.
} 
defined a natural \emph{bundling transformation}
\[
\defm\eta\: \MU_G(X) \lt  \MU_G(EG \x X)
\isoto \MU\big((EG \x X)/G\big),
\]
where $(EG \x X)/G \eqc \defm{X_G}$ is the Borel construction.
It is not hard to see there is a homotopy equivalence
	$EG_+ \wedge_G MU^G_k \hmt BG \x MU_k$
	induced by the map classifying the (nonequivariant) vector bundle
	$(EG \x EU_k^G)/G \to (EG \times BU_k^G)/G$,
and $\h$
takes a class represented by $S^V \wedge X \lt MU^G_{n+|V|}$ to
the class of
\eqn{
	EG_+\mn \us G \wedge\, (S^V  \!\wedge  X) \ \lt\  EG_+\mn \us G \wedge MU^G_{n+|V|} 
	\ \lt \ MU_{n+|V|}
	\mathrlap.
}
When $X = {*}$, the map $\MUG \lt \MU(BG)$ is known to be 
equivalent to the completion map $MU^*_G \lt (MU^*_G)\compl$
with respect to the ideal $\ker(\MUG \to MU^*)$
due to work of Greenlees--May and La Vecchia%
~\cite{greenleesmaycompletion,lavecchiacompletion}.	
When $G = T$ is a torus, 
the composite map $\WUT \lt \MU(BT)$ is often referred to as the 
\emph{universal toric genus}~\cite{buchstaberpanovray}.

For $h \in\{H,K,MU\}$ the \AHSS converging to $h^* BU(n)$
collapses, so that $h^* BU(n)$ is the power series ring over $h_*$
on classes 
$\defm{c_i^{h}} \in \smash{h^{2i}BU(n)}$
called the \defd{Conner--Floyd Chern classes}.
Applying these classes 
to the stable normal bundle of a class $[M \to X] \in mu^G_*(X)$
judiciously and following with the relevant Gysin map
(integration along the fiber for $\HG(X)$, 
the Atiyah--Singer index map for $\KG(X)$), 
one compiles them into ``characteristic number'' maps
from $mu^*_G(X)$, factoring through $MU_*^G(X)$ 
and taking values in 
$\Hom(\K BU, \KG X \big) \iso \KG(X)[[a_1,a_2,a_3,\ldots]]
\eqc \KG(X)[[\defm{\vec a}]]$
or the power series ring 
\[
	\Hom_{h^*}(h^* BU, h^* X_G)\, \iso \, 
	h^*(X_G)  \us {\phantom{{}^*}h_*}{\wh\ox} \, h_* BU \, \iso\,  
	h^*(X_G)[[\vec a]]\mathrlap.
\]
These are also called ``Boardman maps.''\footnote{\ 
	The oldest version, due to Boardman, 
	is $MO_*(-) \lt \Hom\big(\H(BO;\F_2),\H(-;\F_2)\big)$.
}
Most of these maps factor through the bundling map to $\MU(X_G)$
using $\MU(-) \lt \K(-)$ or $\MU(-) \lt \H(-)$.

tom Dieck showed that $\MUG(X) \lt \KG(X)[[\vec a]]$
is injective when $G$ is topologically cyclic and $X$ is a point
or the unit sphere $SV$ in a representation $V$~\cite [Thms.~2\&3]{tomdieckcharacteristicII},
implying $\MUG \lt \MU(BG)$ is injective in these cases as well.
Composing with the map from $\WUG$, this shows manifolds are
determined up to $G$-equivariant cobordism by their K-theoretic 
characteristic numbers,
which Hattori~\cite{hattoriequivariant} proved independently for $G = T$ a torus.
(The case $G = 1$ is due to Stong and Hattori%
~\cite{stongrelations,hattoriintegral}.)
Hattori~\cite[Thm. 1.7]{hattoriequivariant} also proved
for $T$ a torus that $\WUT \lt \MU(BT)$ is monic.
L{\"u} and Wang~\cite{luwangchern}
showed the map $\WUT \lt \H(BT)[[\vec a]]$ is injective as well, 
confirming a conjecture of Guillemin--Ginzburg--Karshon.

Variants of these maps can be defined on $\W^{U:G}_*(X)$
using the equivariant characteristic numbers of the tangent bundle;
their values in $\Hom\big(H_*BU,h^*_G(X)\big)$ 
differ from the composition through $mu^G_*(X)$
by precomposition by the automorphism of $H_*BU$ 
induced by the map $\defm\iota\: BU \lt BU$ 
classifying the ``inverse'' operation 
for the H-space structure on $BU$ classifying the stable Whitney sum.
This corresponds to the exchange 
of the tangent and stable normal bundles.

\subsubsection*{Localization near fixed points}

The fixed point set $M^G$ in a stably complex $G$-manifold $M$
is itself a stably complex $G$-manifold in such a way that its
normal bundle $\nu = \nu_M(M^G)$ carries a natural complex $G$--vector bundle structure.
The compactification $\defm{\bar\nu}$ of $\nu$ 
by its fiberwise visual boundary $\6\nu$
is a stably complex $G$-manifold with boundary
which we tacitly consider 
to be embedded as the closure of a regular neighborhood
of $M^G$ in $M$.
The boundary $\6\nu$ is itself a stably complex $G$-manifold whose
isotropy groups $G_x < G$ for $x \in \6\nu$ are all proper subgroups,
so the pair $(\ol\nu,\6\nu)$ represents an element in a bordism
ring $\defm{\W^{U:G}_*[\ms A\!, \ms P]}$ 
of $G$-manifolds whose isotropy groups
lie in the set $\ms A$ of all closed subgroups of $G$
and whose boundary has isotropy groups in the subset $\ms P$
of \emph{proper} closed subgroups.\footnote{\ 
Such a notion exists for any pair $\ms F \supsetneq \ms G$
of sets of subgroups each stable under subgroup inclusion and conjugacy,
and induction along such pairs is an important proof technique.
There is a tautological long exact sequence
$\cdots \to \Omega^{U:G}_n[\ms G] \to \Omega^{U:G}_n[\ms F] \to 
\Omega^{U:G}_n[\ms F\mn,\ms G] \os\6\to \Omega^{U:G}_{n-1}[\ms G]
\to \cdots$
due already to Conner and Floyd~\cite{connerfloyd}
in the unoriented case
and inducing other such sequences by stabilization, 
localization, and completion.}
The pairs $(M,\0)$ and $(\nu,\6\nu)$
represent the same class in this ring,
and similar reasoning shows every class
is represented by a disc bundle over a $G$--fixed set and its boundary. 
This map is known to be injective for $G$ compact abelian%
~\cite[p.~173]{hamrickossa}.
More generally
there is a ring map
$\defm{\rho_\W}\:\WUG(X) \lt \WUGAP(X)$ natural in $G$-spaces $X$
and similarly there are maps $mu_*^G(X) \lt {mu^{G}_*[\ms A\!, \ms P]}(X)$
and, taking colimits, 
$\defm{\rho_{MU}}\:MU_*^G(X) \lt {MU^{G}_*[\ms A\!, \ms P]}$.

There is a natural additive description for $\WUGAP$.
Let $\defm J$ be a set of
nontrivial irreducible complex $G$-representations,
containing each exactly once up to isomorphism.
Then the interior $E$ of a complex disk bundle $\ol E \to E^G$ 
isotypically decomposes
over each component $N$ of $E^G$ as a direct sum over $V \in J$
of vector bundles of the form $F_V \ox V$.
Such bundles $F_V$ are classified by maps $N \lt BU(m)$,
so that $\W^{U:G}_{n}[\ms A\!,\ms P]$ is isomorphic to the direct sum 
of the groups $MU_{2m_0}\big(\prod_V BU(m_V)\big)$
over lists $(m_0,m_V)_{V \in J}$ of nonnegative numbers 
with $m_0+\sum m_V|V| = m$~\cite[\SS4]{tomdieckintegrality}.

It is also possible to describe $\MUGAP$ in terms of fixed points,
restricting a representative 
$\defm{S^{W+n}} \ceq S^W \wedge S^{n} \to MU^G_{|W|}$ of a class in $MU^G_{n}$
to $S^{|W^G|+n} \to (MU^G_{|W|})^G$.
Since maps to $BU_{|W|}^G$ classify complex $G$-bundles,
and each $|W|$-dimensional $G$-representation can be written
uniquely up to isomorphism as
$\C^{m_0} \+ \Direct_{V \in J} V^{\oplus m_V}$
for some list $(m_0,m_V)_{V \in J}$ with 
$m_0 + \sum_V m_V |V| = |W|$,
one finds 
$\smash{(MU_{|W|}^G)^G}$ is the wedge sum over such lists
of $\smash{MU(m_0) \wedge \big(\prod_{V \in J} BU(m_V)\big)_{\mn +}}$,
and so the fixed-point map represents an element of
the $\smash{\big(|W^G|+n\big){\mn}\th}$
homotopy group of the wedge.
Stabilizing,
let us agree to write $\defm{BU^{\oplus J}}$ 
for the subspace of the product $BU^J$ 
containing those points all but finitely many of whose
components are the basepoint.
It can then be shown, after reshuffling, 
that the geometric fixed point spectrum
$\defm{\Phi^G MU^G}$ decomposes as a wedge, 
indexed by elements
$W$ of the augmentation ideal $\defm{I(G)} \lhd R(G)$ of virtual representations of dimension $0$,
of spectra $\smash{S^{|W^G|} \wedge MU \wedge BU^{\oplus J}}$,
and that $\MUGAP \iso \pi_*(\Phi^G MU^G)$%
~\cite[\SS 2]{tomdieckintegrality}\cite[\SS 4]{sinha2001computations}.


To relate this back to $\MUG$ we must introduce Euler classes.
	Each $G$-representation $V$ is a fiber of the universal $G$-bundle
	the Thom construction converts the fiber inclusion 
	to a $G$-map $\smash{\defm{u_V}\: S^V \longinc MU^G_{|V|}}$.
	The $MU^G$-homological Thom isomorphism
	$\smash{\wt{MU}^G_{*+|V|}(S^V)} \isoto \MUG$
	is given by applying homotopy groups to 
the composition of ${\id} \wedge {u_V}\: MU^G_n \wedge S^V\to 
	\smash{MU^G_n \wedge MU^G_{|V|}}$
and the spectrum multiplication,
which is to say ${u_V}$ represents the Thom class 
of the bundle $V \to *$ in $\MUG$.
The element of $\smash{\pi_{-|V|}(MU^G)} = \smash{MU^G_{-|V|}({*})}$ 
	given by restricting $u_V$ to $S^0 = \{0,\infty\}$ 
	is the Euler class $\defm{e_V}$ of this bundle.
	Smashing, we see $u_V \. u_W = u_{V \oplus W}$
	and $e_V \. e_W = e_{V\oplus W}$.
	
	It is not obvious \emph{a priori} 
	that the $e_V$ are interesting.
	In fact they vanish
	when $V$ has nontrivial invariant subspace,
	for the representative above factors
	through 
	$\smash{S^0 \longinc S^{V^G}}$,
	which is equivariantly nullhomotopic for $|V^G| > 0$.
	Otherwise, however, they are nonzero.\footnote{\ 
		And since they 
		are negative-dimensional,
		they do not come from $mu^G_*$,
		showing that $mu^G_* \lt MU^G_*$ is not surjective
		before stabilization
		and the stabilization maps $mu_*^G(X) \lt mu_{*+|W|}^G(S^W \wedge X)$
		are not all surjective.
	}
	In fact~\cite[Lem.~2.2]{tomdieckintegrality}%
	\cite[\SS4]{sinha2001computations},
	 carefully following through the identification of $\MUGAP$
	with $\pi_*(\Phi^G MU^G)$
	yields graded $MU_*$-module isomorphisms
\[ 
	\pi_* (\Phi^G MU^G)
	= MU_*\Big(\bigvee_{W \in I(G)} S^{|W^G|} \wedge BU^{\oIpG}\Big)
	\iso
	\Z[e_V^{\pm 1}]_{V \in J} \ox MU_*\big(BU^{\oIpG}\big)
\]
such that the image of $e_V \in \MUG$
under $\MUG \to \MUGAP \to \Z[e_V^{\pm 1}] \ox MU_*(BU^{\oplus J})$
is $e_V \ox 1$.
The wedge summand of $\Phi^G MU^G$ 
indexed by a formal difference of representations $W = W' - W'' \in I(G)$
corresponds to the summand 
$\Z e_{W'/(W')^G}\. \smash{e_{W''/(W'')^G}\-} \ox MU_*(BU^{\oIpG})$.
 The ring $MU_*(BU^{\oplus J})$ is 
 $MU_*(BU)^{\ox J}$
 by the K\"unneth formula and the fact $MU_*$ commutes with colimits, 
 which in turn is isomorphic to 
 $MU_*[X_2,X_4,X_6,\ldots]^{\otimes J}$ 
 by the collapse of the Atiyah--Hirzebruch
 spectral sequence for $MU_*(BU)$.
 A usual choice of such $X_{2d} \in MU_{2d}(BU)$ 
 is represented by the
 standard compositions $\defm{\chi_{2d}}\:\CP^d \inc \CP^\infty = BU(1) \inc BU$,
 corresponding in $\MUGAP$ to the disc bundles of the tautological 
 line bundles $\defm{\g^d} \to \CP^d$.
 Thus we finally have an isomorphism
 \[
 \phantom{( V \in J,\  d \geq 1.) \qquad}
 \defm{\chi_{MU}}\:\MUGAP \isoto MU_*[X_{V,\,2d},e_V^{\pm 1}]
 \qquad (V \in J,\  d \geq 1),
 \]
 where $X_{V,\,2d} \in MU_*(BU^{\oIpG})$ 
 comes from the copy of $X_{2d}$ in the $V\th$ $MU_*(BU)$ factor.

The map $\WUGAP \isoto 
	\Direct MU_{2m_0}\big(\prod_V BU(m_V)\big)$
	also extends to a map 
	\[
	\defm{\chi_\W}\:\WUGAP \lt \Z[e_V^{\pm 1}] \ox MU_*(BU^{\oplus J})
	\]
	in a natural way.
	Extending the homotopy class of a map
	 $f\:N \lt \Direct_V BU(m_V)$
	 via the natural map $i\: \Direct_V BU(m_V) \lt BU^{\oplus J}$
	 preserves the stable isomorphism classes of the 
	 isotopic components $E_V \ox V$ of the corresponding disc
	 bundle $E \to N$, but loses the dimension data 
	 $|E| = |N| + \sum_V m_V  |V|$.
	 To remember the dimensions
	 and make the composite
	 $\chi_W \o \rho_\W\:\WUG \to \WUGAP \to \Z[e_V^{\pm 1}] \ox MU_*(BU^{\oplus J})$
	 preserve the grading,
	 we instead take $\chi_\W[f] \ceq \prod_V e_V^{-m_V} \otimes {[i\o f]}$.

	 One might suspect 
	 the composites $\chi_\W \o \rho_\W$
	 and $\chi_{MU} \o \rho_{MU} \o \Psi\: \WUG
	 \lt MU_*[X_{V,\,2d},e_V^{\pm 1}]$ should agree.
	 Instead~\cite[\SS4]{tomdieckintegrality}, 
	 they differ by postcomposition by 
	 the involution $\defm\tau \ceq MU_*(\iota)^{\oplus J} \ox \id$,
	 where $\iota$ is the H-space inverse on $BU$. 
	 To see this, embedding 
	 a stably complex $G$-manifold $M$ in some large representation~$W$,
	 the normal bundles $\nu_M(M^G)$ and $\nu_M(W)|_{M^G}$
	are complementary in $\nu_W(M^G) = \nu_{W^G}(M^G) \+ \nu_W(W^G)|_{M^G}$.
	As $\nu_{W^G}(M^G)$ has trivial $G$-action
	and $\nu_W(W^G)$ is a product bundle,
	the classifying maps of the $V$-isotypic 
	components of the two normal bundles
	 $\nu_M(M^G)$ and $\nu_M(W)|_{M^G}$
	are stably inverse with respect to Whitney sum
	for each $V$.

Further comparison is facilitated by the 
identification of the image of $\WUGAP \to \MUGAP \simto MU_*[X_{V,\,2d},e_V^{\pm 1}]$.
For $d \geq \frac 1 2|V|$
write $\defm{y_{V,\,2d}} \in \W_{2d}^G[\ms A\!,\ms P]$
for the class represented by the disc bundle of 
$\g^{d-|V|/2} \ox V \to \CP^{d-|V|/2}$,
	classified 
	by the inclusion $\CP^{d-|V|/2} \inc \CPi = BU(1)$ 
	in the $V\th$ coordinate.
We write $\defm{Y_{V,\,2d}} = \chi_\W(y_{V,\,2d}) = X_{V,\,2d-|V|}e_V\-$.\footnote{\
	We prefer to index generators by dimension;
	Hanke and Sinha have $Y_{d,V}$ for these same elements.
}
Evidently we may exchange generators
in $MU_*[X_{V,\,2d},e_V^{\pm 1}] = MU_*[Y_{V,\,2d},e_V^{\pm 1}]$,
and then Hanke~\cite[Prop.~4]{hankebordism} proves
that $\chi_\W$ is an isomorphism onto the subring $MU_*[V_{V,\,2d},e_V\-]$.
Thus $\chi_\W$ can be seen as an algebraic localization 
inverting the $e_V\-$.
Hanke~\cite{hankebordism} shows that when $G = T$ is a torus,
then the following is an injective pullback square:
	\quation{\label{eq:Hankesquare}
	\begin{aligned}
		\xymatrix{
			\Omega^{U:T}_* 					\ar[d]_\Psi	
													\ar[r]^(.45){\rho_\W}			&	
			\WUTAP \ar[d]_{\Psi[\ms A\!,\ms P]}		\ar[r]^(.475){\chi_\W}_(.45)\sim	&
			MU_*[V_{V,\,2d},e_V\-] 								\ar@{^{(}->}[d] 			\\
			MU^T_*	 													\ar[r]_(.45){\rho_{MU}}			&
			\MUTAP														\ar[r]_(.475){\tau \o \chi_{MU}}^(.45)\sim&	
			MU_*[V_{V,\,2d},e_V^{\pm 1}] \mathrlap.
		}
	\end{aligned}
}
That is,
$\WUT$ is the intersection of $\MUT$ and $MU_*[Z_{V,\,2d},e_V^{-1}]$
in $MU_*[Z_{V,\,2d},e_V^{\pm 1}]$.

\subsubsection*{Algebraic localization}

	Localization 
	in the geometric sense of $\MUG \lt \pi_*(\Phi^G MU^G)$
	turns out to agree with the 
	algebraic localization inverting the Euler classes $e_V$.
	Let us write $\defm S <  MU^G_*$ 
	for the multiplicative submonoid generated by the Euler classes $e_V$.
	As these become invertible
	in $\MUGAP \iso MU_*[Y_{V,\,2d},e_V^{\pm 1}]$,
	there is an induced homomorphism 
	$S\-MU^G_* \lt \MUGAP$,
	which tom Dieck shows is an 
	isomorphism~{\cite[Thm.~3.1]{tomdieckintegrality}}. 
		Sinha~\cite[Thm.~5.1]{sinha2001computations}
		showed 
		that when $G$ is abelian,
		multiplication by $e_V$ is injective if and only if
		$G$ acts transitively on the unit sphere $SV$,
		which is in particular the case for all representations
		nontrivial on the identity component $T$.
		Thus the localization map $MU_*^T \lt S\-MU_*^T \iso MU_*^T[\ms A\!,\ms P]$
		is injective, but for all nontoral $G$
		the corresponding map is non-injective.
		
		It {can be shown}~\cite[Prop.~4.14]{sinha2001computations}%
		\cite[p.685]{hankebordism}
		 that 
		$\tau \o \chi_{MU} \o \rho_{MU} \o \Psi$
		takes the class of the projectivized representation 
		$P(\C^d \+ V)$ in $\Omega^{U:G}_{2d+|V|}$ 
		to $Y_{V,\,2d} + e_{V^*}^{-d}$ if $|V| = 2$,
		where $V^*$ is the dual representation,
		and to $Y_{V,\,2d}$ if $|V| \geq 4$.
		Write $\defm{Z_{V,\,2d}}$ for this class in either event.
		Then when $G = T$ is a torus, 
		using $\tau \o \chi_{MU} \o \rho_{MU}$ to
		identify $\MUT$ as a subring of
		$MU_*[Z_{V,\,2d},e_V^{\pm 1}] = 	MU_*[Y_{V,\,2d},e_V^{\pm 1}]$,
		Sinha~\cite[Thm.~1.2]{sinha2001computations}
		found inclusions
		\[
		MU_*[Z_{V,\,2d},e_V] \leq 
		MU_*^T \leq 
		S\-MU_*^T \iso 
		MU_*[Z_{V,\,2d},e_V^{\pm 1}]
		\mathrlap.
		\]

The bundling and characteristic number maps
induce $\MUG$-module structures on their codomains,
so one can localize them with respect to $S$ as well.
\quation{\label{eq:squares}
	\begin{aligned}
		\xymatrix@C=1em@R=4em{
			MU^G_*	 	\ar[d]	\ar[r]^\h		&	
			MU^*(BG)		\ar[d]	\ar[r]			&
			H^*(BG)[[\vec a]]	\ar[d]					&
															&
			MU^*(BG)		\ar[d]	\ar[r] &
			H^*(BG)[[\vec a]]		\ar[d]				
			\\
			S\-\MUG			\ar[r]		&	
			S\- MU^*(BG)		\ar[r]		&
			S\-H^*(BG)[[\vec a]]			&
														&
		S\-MU^*(BG)			\ar[r] &
		R(G)[[\vec a]]			
		}
	\end{aligned}
}
It is known that some of these squares are pullbacks.
tom Dieck~{\cite[\SS5]{tomdieckintegrality}}
showed $\MUG$ is the pullback of the square 
 with lower-right corner $S\- MU^*(BG)$ for $G \iso \Z/p$,
 and adopted the phrasing that the values of the bundling map are 
 ``integral'' in $S\- MU^*(BG)$, meaning they appear without denominator. 
For $G = T$ a torus,
Hattori~\cite[Thm. 1.4]{hattoriequivariant} showed that
$\WUT$ is the pullback of the square with lower-right corner 
$S\- RT[[\vec a]]$.
For topologically cyclic $G$, tom Dieck~\cite[Thm.~1]{tomdieckcharacteristicII}
showed that the localized Boardman map 
$S\-\MUG \lt S\- R(G)[\vec a]$ is injective
but for other $G$, the codomain is zero~\cite[Lem.~1]{tomdieckcharacteristicII},
putting sharp limits on this sort of pullback identification.

Again for $T$ a torus,
Darby~\cite[Prop.~3.5]{darby2015torus} shows a version of the pullback square
determining the subring $\W^{U:T,\mathrm{fin}}_*$ 
represented by manifolds with finite fixed point sets
as the pullback of a square with lower-right corner $S\-\MUT$.
Combined with Hanke's result this gives
the same ring as the pullback  of a square with lower-right corner $S\-\MU(BT)$.

Although we cannot have pullback squares of the above sort 
when localization is noninjective,
the vertical maps in the right column are frequently known
to be injective as well.
For example,
Hattori~\cite[Thms.~1.3, 1.7]{hattoriequivariant} 
found that
for $G$ topologically cyclic, 
$S\- \WUG \lt S\-R(G)[[\vec t]]$ is injective, 
and for $T$ a torus, $S\-\WUT \lt S\-\MU(BT)$ is injective.

	\subsection{Explicit computations}

	The most relevant explicit computations of the coefficient ring for $G = T$ a torus
	are those 
	of Sinha and Gusein-Zade~\cite{sinha2001computations,GZ:long}. 
	Sinha determined a set of constraining relations on the image
	of $MU^T_*$ in the localization 
	$MU^T_*[e_V\-] \iso MU_*[Z_{2d,V},e_V^{\pm 1}]$,
	but the tools available did not show this set of relations 
	to be complete unless $T=S^1$.%
	
	Restricting to semifree circle actions, 
	Sinha found explicit presentations%
	~\cite[Thms.~3.6, 3.10]{sinha2005semifree}.
	One defines $\defm\WSF$
	as the bordism ring of semifree stably complex $S^1$-manifolds 
	and finds the natural map to the unrestricted $\smash{\Omega_{\smash*}^{U:S^1}}$
	is injective~\cite[Rmk.~2.5]{sinha2005semifree},
	so we will avail ourselves of our earlier notation. 
	A spectrum $\defm{MU^{\mr{SF}}}$ is 
	defined as $MU^G$ was,
	but using only the trivial representation $\C$,
	the standard representation $t$, and its conjugate $\bar t$.
	Sinha finds $\WSF$ is a free $MU_*$-module
	injecting in $\WSF\AP \iso 
		MU_*[Z_{V,2d},e_V\-]_{V \in \{t,\bar t\},\, d\geq 1}$,
	and similarly $\MUSF$ 
	contains $MU_*[Z_{V,2d},e_V]$
	and injects in 
	 $\smash\MUSF\AP \iso MU_*[Z_{V,2d},\smash{e_V^{\pm 1}}]$, 
	 forming an injective pullback square~\cite[\SS2]{sinha2005semifree}.
	 He finds explicit algebra generators for~$\WSF$
	 and~$\MUSF$ and $MU_*$-module bases for each,
	 which enable him to show a certain 
	 set of geometrically defined relations is complete in each case.
	 	 
	 Returning to the unrestricted case,
		Musin~\cite{musingenerators} found a natural set of generators
		for $\smash{\Omega^{U:S^1}_*}$
		and the fixed point--free ring $\smash{\Omega^{U:S^1}_*}[\ms P]$.
	Gusein-Zade~\cite{GZ:long}
	had earlier found an equational description of $\Omega^{U:S^1}$
	in terms of its image under the injection 
	$\smash{\Omega^{U:S^1}} \lt \smash{\Omega^{U:S^1}}[\ms A\!,\ms P]$
	but this description is stated via equations in a sextuply-indexed 
	array of power series dependent on several previous levels
	of power series, 
	and so cannot be easily applied in practice~\cite[Intro.]{musingenerators}
	to determine if a given equivariant disc bundle actually 
	arises as the normal bundle to the fixed point
	set of a manifold.

		Similarly, in the finite realm,
		Mi{\v{s}}{\v{c}}enko~\cite{mishchenkobordism}
		equationally determined the image of
		$\W^{U:\Z/p}_{\smash{*}}$ under the embedding in $\W^{U:\Z/p}_{\smash{*}}[\ms A\!,\ms P]$, 
		Kosniowski~\cite{kosniowskigenerators} 
		supplied an explicit set of geometric generators,
		and Jack Carlisle has recently found presentations for $\W^{U:\Z/p}_{\smash{*}}$.
		As for homotopical bordism,
		Kriz~\cite{krizZp} found an expression for $\smash{MU_*^{\Z/p}}$
		as a pullback of a \emph{non}injective square
		of maps involving a localization of 
		a quotient of a power series ring,
		Strickland~\cite{stricklandcobordism} 
		then found a presentation for $\smash{MU_{\smash{*}}^{\Z/2}}$,
		Abram--Kriz~\cite{abramkriz} 
		found an algebraic expression for $MU_*^A$
		for $A$ finite abelian,
		and
		Hu has found presentations for $\smash{MU^{\Z/p^r}_{\smash{*}}}$.
	Leaving abelian groups behind, Hu--Kriz--Lu~\cite{hukrizlu} 
	have also computed
	$\smash{MU_*^{\Sigma_3}}$, where $\Sigma_3$ is the symmetric group on three letters.

\section{Semifree actions}
	We have seen characterizations of bordism rings 
	frequently employ homomorphic embeddings into larger rings.
	We adopt a variant of this approach
	to characterize semifree bordism with isolated fixed points.
	
	Given a stably complex $G$-manifold $M$, for each finite list of
	numbers $I$ there corresponds an equivariant Chern class $c^G_I(M)$
	in the Borel cohomology $H^*(EG \ox_G M)$, whose pushforward under the Gysin map
	corresponding to $M \to {*}$ is a so-called Borel \defd{equivariant Chern number}
	$\defm{h_I}(M) = \int_M c^G_I(M) \in H^{* - \dim M} BG$.\footnote{\ 
		These can be compiled into a ring homomorphism 
		$\Omega^{U:T}_* \lt \H BT \,\wh\otimes\, H_* BU$ if one likes,
		but we will not need this point of view in the present note.
	}
	Not just any list $(h_I)$ of elements of $\H BG$ arises as 
	$\big(h_I(M)\big)$ for some $M$, of course;
	for example, when the total degree $\defm{|I|}$ of $c_I$ is less than the real dimension of $M$,
	then $h_I$ is zero because $H^{< 0} BT$ is.

	When $G$ is a torus, 
	the Atiyah--Bott/Berline--Vergne localization formula~\cite{BV1982,AB1984} 
	expresses the Chern numbers 
	$h_I(M)$ entirely in terms of the image of $M$ in $\bundlering$.
	If $\nu\: N \to M^T$ is the normal bundle, viewed as a regular
	neighborhood in $M$,
	then
	\quation{\label{eq:ABBVChern}
		\int_M c^G_I(M) = \int_{M^T}\frac{c^G_I(M)|_{M^T}}{e^T(\nu)}
	}
	as elements of $\H BT$, 
	where $\defm{e^T}(\nu)$ 
	is the Euler class of the 
	induced bundle $\nu \ox \id\: N \ox_T ET \lt M^T \ox_T ET$.
	The expression on the right-hand side of \eqref{eq:ABBVChern},
	 which we will call 
	$\defm{\ell_I}(\nu)$, is defined
	independently of whether $[\nu]$
	lies in the image of the map from $\mfring$ or not,
	but if it does lie in that image, 
	then we know from  \eqref{eq:ABBVChern} that 
	$\ell_I(\nu) = h_I(M) = 0 $ for $|I| < |N|$

	One can ask if these restrictions alone determine the image of $\mfring \lt \bundlering$; 
	that is, given a complex $T$-equivariant bundle 
	$\nu\: N \to X$
	such that $\ell_I(\nu) = 0$ for $|I| < \dim N$,
	does it always ``close up'' to some closed, stably complex $T$-manifold $M$?
	We show the answer is yes
	for certain classes of actions
	and use this fact to recover the associated bordism rings,
	and particularly Sinha's \Cref{thm:semifree bordism}.

	\begin{definition}\label{def:semifree data}
		We write $\defm t$ for the standard one-dimensional complex representation of $S^1 = \U(1)$,
		and $\defm{\bar t}$ for its conjugate,
		and $\defm{V_j} = t^{\oplus n-j} \+ \bar t^{\oplus j}$. 
		We will write $\defm u$ for the image of $t$ under the standard isomorphism
		$\Hom(S^1,S^1) \isoto H^2(BS^1;\Z)$
		taking an irreducible representation $V$ to the first Chern class
		of the associated complex line bundle $ES^1 \ox_{S^1} V \to BS^1$,
		so that the image of $\bar t$ is $-u$.
	\end{definition}
	
	\begin{definition}\label{def:AID}
		\defd{Abstract isotropy data} 
		comprises a family of pairs
		$(V_p,\sp)$, indexed by a finite set $P$,
		of signs $\defm \sp \in \{\pm 1\}$,
		and complex $S^1$-representations $\defm{V_p}$
		all of one common dimension.
		We say abstract isotropy data $(V_p,\sp)_{p \in P}$
		is \defd{semifree} if each $V_p$ is isomorphic to 
		one of the representations $V_j$ ($1 \leq j \leq n$)
		for some fixed natural number $\defm n > 1$.
		In this case we define $\defm q\: P \lt \{0,1,\ldots,n\}$
		by $V_p \iso V_{q(p)}$.
		The \defd{isotropy data} of an oriented stably complex $S^1$-manifold $M$
		with isolated fixed points is
		the \AID $(T_p M,\s_p){}_{p \in \smash{M^{S^1}}}$,
		where $\s_p$ is $1$ if the given orientation of $M$
		agrees with the orientation on $T_p M$ induced by the stable complex
		structure and $-1$ otherwise.
		
	\end{definition}
	
	If \AID are viewed as multisets by forgetting $P$
	but remembering multiplicity, 
	they form a semiring with addition given by disjoint 
	union and multiplication given by 
	$(V,\s)\defm{\,\.\,}(V',\s') = (V \oplus V',\s\s')$.
	This semiring is generated by $(V,1)$ and $(V,-1)$,
	where $V$ runs over (isomorphism classes of) 
	irreducible nontrivial representations.
	The ring localization of this semiring
	can be viewed as the ring of ``\AID up to bordism,'' 
	defined by quotienting by $-(V,\s) \sim (V,-\s)$
	and writing $\defm{\s V}$ for the corresponding 
	class.
	This ring is a polynomial $\Z$-algebra in the 
	nonzero irreducible representations.
	Restricting to semifree \AID,
	the subsemiring is generated by $(t,\pm 1),(\bar t, \pm 1)$
	and the ring is $\Z[t,\bar t]$.
	%
	%
	%
	%
	%
	We will do a bit more than characterize the	semifree 
	bordism ring, in fact computing the subsemiring
	of geometrically realized isotropy data within the \AID.
	
	\medskip
	
	Given a compact, oriented, stably complex, semifree $S^1$-manifold $M^{2n}$
	with isolated fixed points and isotropy data $(V_p,\sp)_{p \in P}$,
	the formula \eqref{eq:ABBVChern} specializes to 
	the identities
	\[
	0 
	= 
	\int_M \mn c_i(T M) 
	=\!
	\sum_{p \in M^{S^1}} \frac
	{c_i(\overbrace{u,\ldots,u}^{n-q(p)},\overbrace{-u,\ldots,-u}^{q(p)})}
	{\sp \,\mn u^{n-q(p)}(-u)^{q(p)}}
	=\!  
	\sum_{p \in M^{S^1}} \frac{ c_i(\overbrace{1,\ldots,1}^{n-q(p)},\overbrace{-1,\ldots,-1}^{q(p)})
	}
	{\sp (-1)^{q(p)}u^{n-i} }
	\ \ \ \ \ (0 \leq i \leq n-1),
	\]
	where $\defm{c_i}$ denotes both the $i\th$ Chern class
	and the $i\th$ elementary symmetric polynomial
	since the symbol $\s$ is taken.
	Evidently we can multiply the $i\th$ identity through by $u^{n-i}$
	without changing its content,
	so these identities are really statements 
	about the integers $\sp$ and $q(p)$.
	
	\begin{definition}
		The \defd{ABBV identities} for semifree abstract isotropy data
		are the equations 
		\quation{\label{eq:ABBV semifree}
			\phantom{		\qquad	(0 \leq i \leq n-1).}
			\defm{I_i} = \sum_{p \in P} \sp (-1)^{q(p)}
			c_i(\overbrace{1,\ldots,1}^{n-q(p)},\overbrace{-1,\ldots,-1}^{q(p)})
			=
			0
			\qquad	(0 \leq i \leq n-1).
		}
	\end{definition}

	We will demonstrate the converse, that if \AID
	satisfies the ABBV identities, 
	then it is the isotropy data of some action,
	by showing that the identities prescribe the function $q$ 
	so rigidly that all possible $q$ arise from disjoint unions 
	of the following known examples.

	\bex\label{eg:S2}
	Endow $S^2$ with the standard action and the complex structure
	of $\CP^1$,
	so that the isotropy representation at the north pole $z$ is $t$
	and that at the south pole $-z$ is $\bar t$.
	Let $\e_1,\ldots,\e_n$ each be $1$ or $-1$.
	Endowing the direct power $(S^2)^n$ with the product complex action
	and diagonal $S^1$-action,
	the isotropy representation at the fixed point $(\e_1 z,\ldots, \e_n z)$ 
	is $V_j$, where $j = \smash{\big|\{k : \e_k = -1\}\big|}$,
	so that there are precisely
	$\smash{n \choose j}$ fixed points $p$
	such that $\smash{T_p (S^2)^n} \iso V_j$.
	One has $\s_p = 1$ at each fixed point $p$.
	\eex
	
	%
	%
	%
	
	\bex\label{eg:repsphere}
	The real $S^1$-representation $V_j \+ \R$
	carries a natural stable complex structure given by restriction of 
	the trivial bundle ${\ul V_j} \+ \ul\R \+ \ul\R = {\ul V_j} \+ \ul \C$,
	inducing a stably complex $S^1$-manifold structure on the unit sphere $S(V_j \+ \R)$.
	The $S^1$-action on  $S(V_j \+ \R)$ has two fixed points, one with isotropy data 
	$(V_j,1)$ and the other with $(V_j,-1)$.  
	\eex
	
\nd	Of course only the sphere-powers of \Cref{eg:S2} are necessary 
	to generate the bordism ring, the manifolds $S(V_j \+ \R)$
	being nullbordant by definition, but both will be necessary to 
	generate the semiring of abstract data satisfying
	the ABBV identities \eqref{eq:ABBV semifree}.
	
	Let us now derive a more transparent form of these identities.
	In all of the following, $n$ is fixed.
	Recall that for natural numbers $j,k$, the expression
	\[
	{j \choose k} = \frac {\, \,\,j (j-1) \cdots (j-k+1)}{k (k-1) \cdots
		1 \phantom{j-k+{})}}
	\]
	makes sense even if $k \geq j$,
	yielding the empty product $1$ in case of equality and $0$ if $k > j$.
	Thus we have a polynomial equation
	\[
	\defm{C_i}(j) \ceq 
	c_i(\overbrace{1,\ldots,1}^{n-j},\overbrace{-1,\ldots,-1}^{j})
	=
	\sum_{k = 0}^i{n-j \choose i-k}(1)^{i-k}{j \choose k} (-1)^k 
	=
	\sum_{k = 0}^i \defm{a_{i,k}} j^i \in \Q[j]
	\]
	for some rational numbers ${a_{i,k}}$ with $a_{i,i} \neq 0$.
	Thus we can write
	\begin{equation}\label{eq:IS}
	I_i = \sum_{p \in P} \s_p(-1)^{q(p)} C_i\big(q(p)\big)
		= \sum_{k=0}^i a_{i,k} \underbrace{\sum_{p \in P} \s_p(-1)^{q(p)} q(p)^k}_{\defm{S_k}}\mathrlap.
	\end{equation}
	Assuming the ABBV identities $I_i = 0$ all hold,
	since the $a_{i,i}$ are nonzero,
	one finds inductively that $S_i = 0$ for each $i$ as well.

%
	%
	%
	%
	
	Now we define notation for the multiplicities of the different representations: 
	\eqn{
		\defm{m_j^{+\phantom{-}}}\!\!\!\! &\ceq \big|\{p \in P : 
		\sp = \phantom{-}\!1
		\mbox{ and } q(p) = j\}\big|,
		\\
		\defm{m_j^{-\phantom{+}}}\!\!\!\! &\ceq \big|\{p \in P : 
		\sp = -\!1
		\mbox{ and } q(p) = j\}\big|,
		\\
		\defm{m_j^{\phantom{-}}} &\ceq m_j^+ - m_j^-
	}
	for $0 \leq j \leq n$. 
	Gathering terms by $q$-value, the identities $S_i = 0$ then become
	\quation{\label{eq:mi}
		\phantom{
			\qquad
			(0 \leq i \leq n-1).}
		\sum_{j = 0}^n (-1)^j m_j j^i = 0
		\qquad
		(0 \leq i \leq n-1)\mathrlap,
	}
which can be written in matrix form as follows:
	\quation{\label{eq:matrix}
		\begin{bmatrix*}[l]
			1 & -1 & 1 & \cdots & (-1)^{n-1} \phantom{i^{-1}} \\
			1 & -2 & 3 & \cdots & (-1)^{n-1} n \phantom{{}^{n-1}}\\
			1 & -4 & 9 & \cdots & (-1)^{n-1} n^2 \phantom{{}^{n}}\\
			\,\vdots & \phantom{-}\,\vdots & \,\vdots & \ddots & 
			\phantom{(-1)^n}\vdots \\
			1 &  -2^{n-1} &3^{n-1} & \cdots & (-1)^{n-1} n^{n-1}
		\end{bmatrix*}
		\begin{bmatrix}
			m_1 \\ m_2 \\ m_3\\ \ \vdots\,\,\, \\ m_n
		\end{bmatrix}
		=
		\begin{bmatrix}
			m_0 \\ 0 \\ 0\\ \vdots \\ 0	
		\end{bmatrix}\!\!.
	}
	The determinant of the square matrix 
	is up to sign a Vandermonde determinant, hence nonzero,
	so the matrix is invertible.
	Multiplying both sides by the inverse shows
	$m_1,\ldots,m_n$ are uniquely determined by $m_0$,
	in fact scalar multiples of $m_0$,
	and	we claim that in fact $m_j = m_0 \smash{n \choose j}$.
But clearly these multiplicities are realized 
	by the union $M$ of 
	$m_0$ disjoint copies of the standard $(S^2)^n$
	of \Cref{eg:repsphere},
	and since the isotropy data of $M$ 
	satisfy \eqref{eq:ABBV semifree}, they also satisfy $\eqref{eq:mi}$.\footnote{\ 
	To see this without topology,
	expand $(1+x)^n$ by the binomial theorem, differentiate $i$ times
	with respect to $x$, and 
	and evaluate at $x =- 1$ to find
	\quation{\label{eq:binomial}
		\phantom{			\qquad
			(0 \leq i \leq n-1).}
		0 = \sum_{j=0}^n {n \choose j} \frac {j!}{(j-i)!} (-1)^{j-i} 
		\qquad
		(0 \leq i \leq n-1)\mathrlap.
	}
	Since $j!/(j-i)!$ is a polynomial of degree $i$ in the variable $j$, 
	the right-hand side of \eqref{eq:binomial} 
	is a $\Q$-linear combination of the quantities
	\[
	\defm{D_i} = \sum_{j=0}^n  (-1)^j {n \choose j} j^i\mathrlap.	
	\]
	That $D_0 = 0$ is \eqref{eq:binomial} for $i = 0$,
	and that the other $D_i$ vanish follows inductively from 
	\eqref{eq:binomial} by subtracting off multiples of $D_k = 0$ for $k < i$.}
	We thus have shown 
	any semifree \AID satifsying \eqref{eq:ABBV semifree}
	satisfies
	\quation{\label{eq:mjchoose}
		m_j = {n \choose j} m_0.\footnote{\ 
			Plugging this back into the expression \eqref{eq:IS},
			one finds the combinatorial identities
		\[
			\phantom{\qquad (0 \leq i \leq n-1)}
			\sum_{j=0}^n \sum_{k=0}^i (-1)^{j+k}	
			{n \choose j} {n-j \choose i-k}{j \choose k}
			= 0
			\qquad (0 \leq i \leq n-1)\mathrlap,
		\]
		which do not seem to be otherwise obvious.
		It would be nice to have a combinatorial proof.
		}
	}
	
	To realize arbitrary \AID satisfying \eqref{eq:mjchoose},
	assume first that $m_0$ is nonnegative,
	and let $M'$ be the disjoint union of $m_0$ copies of $(S^2)^n$
	with the standard orientation,
	so that for each $j$ 
	the isotropy data contains $m_j = m_{\smash j}^+ - m_{\smash j}^-$ instances of $(V_j,1)$ 
	and none of $(V_j,-1)$.
	We will realize the isotropy data successfully
	if we can add ${m_{\smash j}^-}$ instances each of $(V_j,1)$ and $(V_j,-1)$
	for each $j$,
	which we can do by taking the disjoint union of $M'$
	with $m_{\smash j}^-$
	copies of $S(V_j \+ \R)$ from \Cref{eg:repsphere}.
	To instead handle the case $m_0 < 0$, 
	one can simply reverse orientations across the board.
	Thus, up to renaming fixed points,
	$(V_p,\s_p)_{p \in P}$ is the isotropy data of
	\[\label{eq:semifree construction}
	M 
	\,\ceq\, 
	\case{\dsp\coprod_{m_0} \phantom{-}(S^2)^n \,\dis\,\, \coprod_{j = 0}^n \coprod_{\  \smash{m_j^-}\vphantom{0}} S(V_j \+ \R),&
		\mbox{if }m_0 \geq 0,\\
		\dsp\coprod_{m_0} -(S^2)^n \,\dis\,\, \coprod_{j = 0}^n \coprod_{\ \smash{m_j^+}\vphantom{0}} S(V_j \+ \R),&
		\mbox{if }m_0 \leq 0.
	}
	\]
	
	We have proved the following.
	
	\begin{theorem}\label{thm:semifree isolated}
		Any semifree \AID  $(V_p,\sp)_{p \in P}$
		satisfying the ABBV identities \eqref{eq:ABBV semifree}
		is the isotropy data of 
		a compact, oriented, stably complex, semifree $S^1$-manifold $M^{2n}$
		with isolated fixed points.
	\end{theorem}
	
	Now recall from the discussion after \Cref{def:AID}
	that the ring of local data for geometric semifree $S^1$-equivariant
	complex bordism is 
	$\Z[t,\bar t]$.
	The image of $S(V_j \+ \R)$ in this ring is zero,
	and that of $S^2$ is $t + \bar t$.
	Since the map from the geometric bordism ring to $\Z[t,\bar t]$
	is injective, \Cref{thm:semifree isolated} shows its image is $\Z[t+\bar t] = \Z[S^2]$,
	yielding Sinha's \Cref{thm:semifree bordism}.

%
%
%
%

\bs

\begin{remark}
Similar reasoning yields another result of Sinha~\cite[Thm.~1.6]{sinha2001computations},
namely that any stably complex
$4$-dimensional $S^1$-manifold with precisely three fixed points
is equivariantly
cobordant to the projectivization 
$P(\C \+ V \+ W)$ for some irreducible $S^1$-representations $V$ and $W$.
Sinha mostly calculates using Euler classes in $S\-MU^{S^1}_*$, 
but the result also follows by applying 
the ABBV formula~\eqref{eq:ABBVChern}
to $c_0$ and $c_1$
to determine relations amongst the six
weights and three signs $\s_p$
and comparing those for $P(\C \+ V \+ W)$.
The same result also follows from 
	theorems of Jang~\cite[Thm.~7.1 (resp. Thm.~1.1)]{jang2018circle}
	classifying possible isotropy data
	for $S^1$-actions 
	on compact, oriented $4$-manifolds
	with three (resp. finitely many) fixed points.
\end{remark}

	\bs
	
	\nd \emph{Acknowledgments}.
	The impetus for this note is joint work in preparation
	with Elisheva Adina Gamse and Yael Karshon approaching a related question for GKM-actions,
	which in turn arose from a preliminary sketch by Karshon,
	Viktor L. Ginzburg, and Susan Tolman dating to the late 1990s. 
	The author would like to thank Dev Sinha and Bernhard Hanke for 
	their generosity in discussing their work,
	Igor Kriz for updating him on recent work in the case $G$ is finite,
	Jeffrey Lagarias for encouragement
	and for hunting down the \emph{locus classicus} 
	for the failure 
	of equivariant transversality,
	Omar Antol\'in Camarena and Sean Tilson
	for discussing $RO(G)$-grading and possible references,
	and Joanne Quigley for proofreading.

	{\footnotesize\bibliography{bibshort}}
	
	\nd\footnotesize{\textsc{Department of Mathematics, Imperial College London
		}\\
		\url{j.carlson@imperial.ac.uk}
	}
\end{document}